# Distributionally Robust Chance Constrained Programming with Generative Adversarial Networks (GANs)

Shipu Zhao, Fengqi You[*]

College of Engineering, Cornell University, Ithaca, New York 14853, USA

## Abstract

This paper presents a novel deep learning based data-driven optimization method. A novel generative adversarial network (GAN) based data-driven distributionally robust chance constrained programming framework is proposed. GAN is applied to fully extract distributional information from historical data in a nonparametric and unsupervised way without a priori approximation or assumption. Since GAN utilizes deep neural networks, complicated data distributions and modes can be learned, and it can model uncertainty efficiently and accurately. Distributionally robust chance constrained programming takes into consideration ambiguous probability distributions of uncertain parameters. To tackle the computational challenges, sample average approximation method is adopted, and the required data samples are generated by GAN in an end-to-end way through the differentiable networks. The proposed framework is then applied to supply chain optimization under demand uncertainty. The applicability of the proposed approach is illustrated through a county-level case study of a spatially explicit biofuel supply chain in Illinois.

# Introduction

Deep learning utilizes computational models with multiple processing layers to learn data representations with multi-level abstraction.[1-4] Recently, numerous deep learning algorithms have been proposed to solve problems in various areas such as data analytics, image processing and machine translation.[5,6] Leveraging deep learning methods to extract useful information from data to support decision-making has gained increasing popularity.[7] Compared with other machine learning methods, complex correlations and even some hidden modes within data can be found by deep learning methods.[8,9] In data-driven optimization frameworks, uncertainty is usually modeled based on available data.[10-14] Therefore, by efficiently extracting high-level features of data, deep learning has become a powerful method to model uncertainty in the data-driven optimization frameworks.[15,16]

Mathematical programming techniques for optimization under uncertainty have achieved great success in various applications such as process design, scheduling, control and supply chain management.[17-29] Among these techniques, robust optimization, stochastic programming and chance constrained programming are the most popular frameworks. The idea of robust optimization is seeking the optimal solution under the worst case.[30] The worst case is taken with respect to an uncertainty set which models uncertain parameters. The goal of stochastic programming is to find the optimal decision that maximize/minimize the expectation value of the objective function while staying feasible for all the possible scenarios of uncertain parameters.[31] The uncertain parameters are modeled with a probability distribution which is assumed as given a priori in stochastic programming. In chance constrained programming, the objective is to find the optimal solution satisfying a constraint by at least a pre-specified probability.[32,33]

Conventional chance constrained programming often assumes perfect knowledge about the probability distribution of uncertain parameters. However, in practice, such a probability distribution is unknown, and the approximations are usually not fully trusted. Moreover, the optimal solution may be very sensitive to the ambiguous probability distribution, which may lead to suboptimal solutions.[34] To handle ambiguous probability distributions in uncertain parameters, a set of possible probability distributions are taken into consideration instead of assuming only one distribution is correct. Thus, distributionally robust chance constrained programming is proposed to hedge against the distributional ambiguities. The set of possible probability distributions is recognized as ambiguity set. In distributionally robust chance constrained programming problems, the chance constraints are required to be satisfied for each probability distribution in the ambiguity set,[35,36] and tractable reformulations are required for solving this type of problems.[35] One way is to reformulate a distributionally robust chance constrained programming problem into a conventional chance constrained programming problem. In general, chance constrained programming problems are also considered intractable due to its non-convexity and complexity. Sample average approximation (SAA) is an intuitive and popular method for solving chance constrained programming problems.[37] The SAA problem is formulated by replacing the latent random distributions with their empirical counterparts constructed using the drawn samples.[38] Therefore, distributionally robust chance constrained programming problems can hedge against ambiguous distributions and can be solved by SAA with tractable reformulations.

Another data-driven optimization method considering distributional ambiguity is disributionally robust optimization.[39-41] Similar with the distributionally robust chance constrained programming, disributionally robust optimization considers a set of possible probability distributions for uncertain parameters with an ambiguity set. However, distributionally robust optimization usually considers worst-case expectation in the objective function and the worst case is taken with respect to the objective function.[42,43] In other words, uncertain parameters are involved in the objective function. In the distributionally robust chance constrained programming, uncertain parameters are involved with constraints and the worst case is taken with respect to these constraints instead of the objective function.[35]

Unsupervised learning with generative adversarial network (GAN) has proven quite successful recently.[44] GAN has been widely adopted in data generation and sampling.[45,46] The generative model from a learned GAN can easily serve as an estimated density model of the training data.[47] Therefore, GAN is an efficient way to model uncertainty in the data-driven optimization frameworks. The idea for GAN is to train a generator and a discriminator simultaneously, and they will compete during training steps. The competition drives both generator and discriminator to increase their ability until Nash equilibrium is reached. Different from other deep generative models which usually adopt approximation methods and prior assumptions for intractable functions or inference, GAN does not require a priori approximation or assumption. Besides, it can be trained in an end-to-end approach through the differentiable networks. GAN can extract distributional information from training data and generate data samples simultaneously, which fits very well with SAA method in the data-driven distributionally robust chance constrained programming.

To the best of our knowledge, there are very few existing studies integrating deep learning methods with data-driven optimization frameworks and exploring the efficiency of deep learning to model uncertainty in data-driven optimization methods. To be more specific, no existing study has employed GAN to model uncertainty in data-driven distributionally robust chance constrained programming. To fill this knowledge gap, we are faced with two challenges. The first one is how to integrate GAN with distributionally robust chance constrained programming problems in order to leverage power of deep learning for optimization under uncertainty. The second one is how to solve the resulting distributionally robust chance constrained programming problems, which are considered intractable in general. In other words, we need to find an appropriate tractable reformulation which can be integrated with GAN well. Since GAN utilizes deep neural networks, complicated data distributions and modes can be learned. Thus, GAN can model uncertainty in an efficient and accurate way in the optimization under uncertainty frameworks.

In this work, a novel framework is proposed for GAN based data-driven distributionally robust chance constrained programming. GAN is employed to model uncertainty in this data-driven optimization framework from historical data in a nonparametric and unsupervised way. Moreover, by utilizing the power of deep neural networks, GAN can model uncertainty in a more efficient and accurate way in the optimization under uncertainty frameworks. Different from other density estimation methods, GAN does not require a priori approximation or assumption. Distributionally robust chance constrained programming problems are reformulated into conventional chance constrained programming problems and finally solved by the SAA method. The required data samples in SAA method are generated by GAN in an end-to-end way through the differentiable networks. Thus, GAN can extract distributional information from historical data and generate data samples simultaneously. The proposed framework is then applied to supply chain optimization under demand uncertainty to illustrate the applicability. The demand for each customer in each time period is treated as an uncertain parameter with an unknown ambiguous probability distribution. The unknown probability distribution is considered belonging to an ambiguity set which is constructed using historical data. Facility locations, production capacities, external supplies, transportation, procurement and inventory decisions need to be determined to satisfy the demand constraints with at least a pre-defined probability in the worst case. The worst case is taken over all the possible probability distributions in the ambiguity set. The objective is to minimize the total cost including capital cost, operating cost, transportation cost, procurement cost, inventory cost and external product cost. The applicability of the proposed approach is further illustrated through a county-level case study of a spatially explicit biofuel supply chain in Illinois.

The novelties of this work include: (1) a novel GAN based data-driven distributionally robust chance constrained programming framework; (2) a deep learning based framework for supply chain optimization under uncertainty; (3) a novel application to spatially explicit biofuel supply chain optimization in Illinois.

The remainder of this article is organized as follows. We first present a brief introduction to GAN, followed by the framework of GAN based data-driven distributionally robust chance constrained programming. To illustrate the proposed framework, a motivating example is given. Next, we apply the proposed framework to supply chain optimization under uncertainty with problem statement, model formulation and solution algorithm. To further illustrate the applicability of our approach, a case study of a county-level spatially explicit biofuel supply chain in Illinois is presented. Conclusions are given in the last section.

## Generative Adversarial Network

GAN is a recently introduced method for training generative models with neural networks.[47] This approach sidesteps some of the common problems among generative models and adopts a simple training regime. The basic idea of GAN is to train a generator and a discriminator simultaneously, and they will compete towards Nash equilibrium during the training steps. Following the idea, GAN trains a generative model through a minimax game.[48,49] The generative model is pitted against a discriminative model that learns to determine whether a sample is from the real distribution or from a synthetic data distribution. Figure 1 illustrates the main training and working process of GAN. The generative model uses random noise as input and generates synthetic data samples, while the discriminative model is trying to detect the counterfeit samples. Competition in this game drives both generator and discriminator to improve their ability until they reach Nash equilibrium.[44,50,51] The training procedures for GAN aim to find a Nash equilibrium of a non-convex game with continuous, high-dimensional parameters. These properties make training a GAN very difficult, and sometimes the training procedures do not converge. Therefore, a number of training techniques are adopted including minibatch discrimination, historical averaging, normalization and one-side label smoothing.

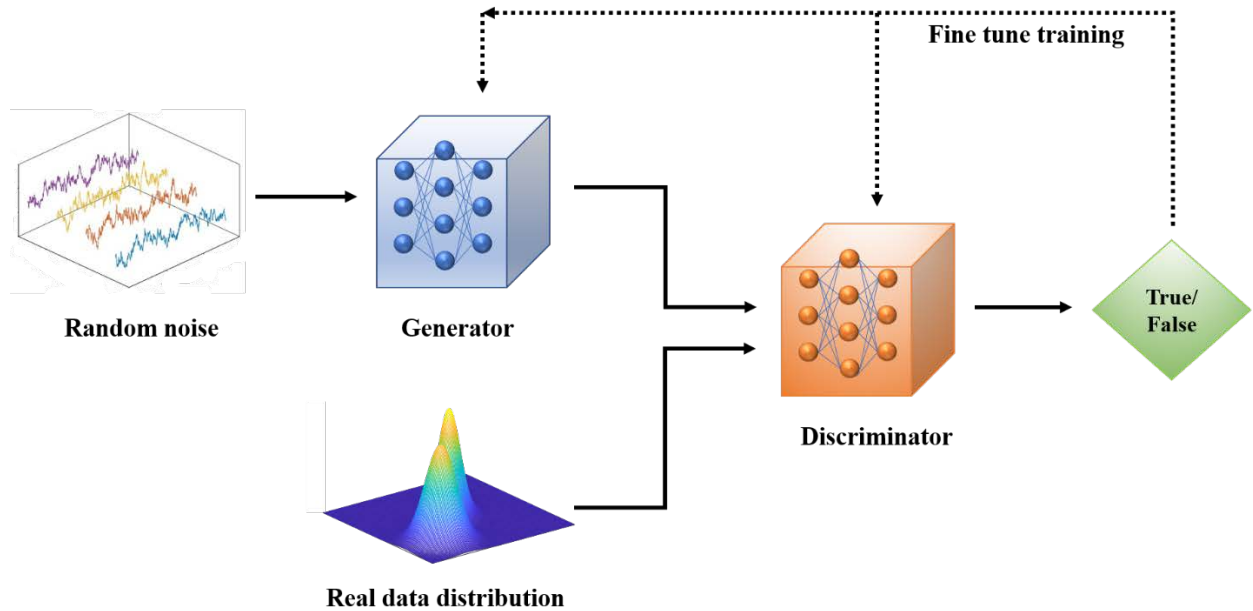

Figure 1   GAN training and working process for generator and discriminator.

Mathematically, the goal for training GAN is to learn a generator distribution $P_G$ which matches the real data distribution $P_{data}$. Instead of trying to give an explicit probability distribution function or assign probability to each data point explicitly, GAN learns a generator network $G$ generating samples $G(z)$ under generator distribution $P_G$ from random noise $z$. Random noise variable $z$ is assumed with a pre-defined probability distribution $P_{noise}(z)$. Thus, the generator network can be regarded as a transformation which transforms a random noise variable to a sample under generator distribution. For each generated sample $x$ from the generator network, the discriminator will determine whether it is from the real data distribution $P_{data}$ or it is a fake one based on the output of the discriminator $D(x)$. Therefore, the learning objective for GAN $V(D,G)$ is given by eq. (1).[47]

$$\min_{G} \max_{D} V(D,G) = E_{x \sim P_{data}(x)} \log D(x) + E_{z \sim P_{noise}(z)} \log\left(1 - D(G(z))\right) \qquad (1)$$

# GAN-based Data-Driven Distributionally Robust Chance Constrained Programming

*Distributionally robust chance constrained programming*

In conventional chance constrained programming, the goal is to find the minimum cost or maximum profit satisfying some constraints (known as chance constraints) by at least a pre-specified probability.[52] Conventional chance constrained programming assumes the probability distribution of uncertainty is known *a priori*. However, it is always challenging to know the exact probability distribution of uncertain parameters in practice, and the optimal solution may be very sensitive to the ambiguous probability distribution. Therefore, distributionally robust chance constrained programming is proposed to address these drawbacks.[35] Instead of considering only one probability distribution, we seek to consider a set of possible probability distributions which are represented as an ambiguity set. The chance constraints should be satisfied by at least a pre-defined probability for all the possible probability distributions in the ambiguity set. In other words, under the worst case, which is taken over the ambiguity set, the chance constraints should be satisfied by at least a pre-defined probability. The distributionally robust chance constrained programming problem is shown as problem (P0).

(**P0**) $\quad \min f(x)$ (2)

$\quad s.t. \quad \inf_{P \in D} \mathbb{P}_{\varsigma \sim P}(C(x,\varsigma) \leq 0) \geq 1 - \alpha$ (3)

$\quad x \in X$ (4)

In eq. (2), $f(x)$ is the objective function to be minimized and $x$ represents the decision variables in a bounded convex set $X$ in constraint (4). Constraint (3) represents the distributionally robust chance constraint where $D$ is the ambiguity set, which can be constructed based on the following subsection. $P$ is the probability distribution for uncertain parameter $\varsigma$. The ambiguity set $D$ usually contains the true probability distribution. Constraint (3) indicates that constraint $C(x, \varsigma) \leq 0$ should be satisfied with at least a pre-defined probability $1-\alpha$ over all the possible probability distributions in the ambiguity set $D$.

*Ambiguity set construction*

To construct the ambiguity set $D$ in the distributionally robust chance constrained programming problem (P0), the concept of distance is used. The distance between two vectors defines how close they are in the vector space. Similarly, we seek to find the distance between probability distributions over the set of all the probability distributions. One commonly adopted way to define the distance between probability distributions is the $\phi$-divergence shown in eq. (5).[53]

$$\phi(f, f_0) = \int_\Omega g\left(\frac{f(\xi)}{f_0(\xi)}\right) f_0(\xi) d\xi \tag{5}$$

where $f$ and $f_0$ denote the true density function and its estimated density function, respectively. To satisfy eq. (5), $g(\cdot)$ should be a convex function on $\mathbb{R}^+$ with the following properties given in eq. (6) – (9).[35]

$$g(1) = 0 \tag{6}$$

$$0 \cdot g(x/0) = x \lim_{p \to +\infty} g(p)/p, \ x > 0 \tag{7}$$

$$0 \cdot g(x/0) = 0, \ x = 0 \tag{8}$$

$$g(x) = +\infty, \ x < 0 \tag{9}$$

One example of such $g(\cdot)$ is given in eq. (10).

$$g(x) = |x-1|, \ x \geq 0 \tag{10}$$

Except for the function given in eq. (10), there are other functions satisfying properties described in eq. (6) – (9). However, the function given in eq. (10) can provide a closed form of the revised risk level,[35] which will be illustrated in the next subsection.

Based on the distance defined by $\phi$-divergence between probability distributions, we can construct the ambiguity set as follows.[54]

$$D = \{P \in M : \phi(f, f_0) \leq d, f = dP/d\xi\} \tag{11}$$

In eq. (11), $M$ is the set of all the probability distributions, $P$ denotes probability distribution and $f$ is the corresponding density function. Thus, for an arbitrary probability distribution in ambiguity set $D$, the $\phi$-divergence between the probability distribution and the empirical probability distribution given by its density $f_0$ should not be greater than a pre-defined value $d$. Parameter $d$ can reflect the range of the ambiguity set and represent the risk-aversion level of decision makers. To be risk-averse, one may choose a large value of $d$, leading to a large ambiguity set.

## *Reformulation of distributionally robust chance constrained programming*

In general, distributionally robust chance constrained programming problem (P0) is considered intractable.[55] Only distributionally robust chance constrained programming problems with specific formations can be reformulated into tractable formations and then solved efficiently. Since the ambiguity set has been defined, we seek to find the reformulation of problem (P0) to get close to the tractable formation. From existing studies, it has been proved that the distributionally robust chance constraint (3) is equivalent to a conventional chance constraint with a modified risk level $\alpha'$ if the ambiguity set is constructed based on $\phi$-divergence, as shown in eq. (12) – (15).[35]

$$\inf_{P \in D} \mathbb{P}_{\varsigma \sim P}(C(x,\zeta) \leq 0) \geq 1-\alpha \Leftrightarrow \mathbb{P}_{\varsigma \sim P_0}(C(x,\zeta) \leq 0) \geq 1-\max\{0,\alpha'\} \tag{12}$$

$$D = \{P \in M : \phi(f, f_0) \leq d, f = dP/d\xi\} \tag{13}$$

$$\alpha' = 1 - \inf_{z, z_0 \in H} \left\{ \frac{g^*(z_0+z) - z_0 - \alpha z + d}{g^*(z_0+z) - g^*(z_0)} \right\} \tag{14}$$

$$H = \{z, z_0 : z > 0, z_0 + \pi z \leq l_g, m^- \leq z_0 + z \leq m^+\} \tag{15}$$

where $\alpha'$ is the modified risk level, $P_0$ is the empirical probability distribution and function $g^*(\cdot)$ is the conjugate of function $g(\cdot)$. $l_g$, $\pi$, $m^-$ and $m^+$ are limit values related to function $g(\cdot)$ and its conjugate $g^*(\cdot)$.

Eq. (12) – (15) are not very straightforward and we cannot employ them directly when solving a distributionally robust chance constrained programming problem. Under special cases of $\phi$-divergence, the modified risk level can be obtained in closed forms.[35] For example, if function $g(\cdot)$ is given as eq. (10), the modified risk level $\alpha'$ can have a closed form shown as eq. (16).

$$\alpha' = \alpha - \frac{d}{2}, \ g(x) = |x-1| \tag{16}$$

Therefore, the distributionally robust chance constrained programming problem (P0) can be reformulated into a conventional chance constrained programming problem (P1),

(**P1**) $\quad \min f(x) \tag{17}$

$\quad\quad$ s.t. $\mathbb{P}_{\varsigma \sim P_0}(C(x,\zeta) \leq 0) \geq 1-\alpha' \tag{18}$

$\quad\quad\quad x \in X \tag{19}$

where the ambiguous probability distribution $P$ is replaced by the empirical probability distribution $P_0$ which can be estimated from available data.

*GAN based empirical probability distribution estimation*

Since we have reformulated problem (P0) to a conventional chance constrained programming problem (P1), the next step is to estimate the empirical probability distribution $P_0$ from available data. In this work, we adopt GAN to estimate the empirical probability distribution. Theoretical results for training GAN are listed: (1) For a given generator $G$ with generator distribution $P_G$, the optimal discriminator is given in eq. (20).[47]

$$D_G^*(x) = \frac{P_{data}(x)}{P_{data}(x) + P_G(x)} \tag{20}$$

(2) The global optimum of the minimax game for training a GAN is achieved if and only if $P_{data} = P_G$ and the optimal $D^* = 1/2$.[47] (3) $P_G$ will converge to $P_{data}$, if $G$ and $D$ have enough capacity, and at each step of training, the discriminator reaches its optimum of the given $G$, and $P_G$ is updated based on the training objective given in eq. (1).[47]

From these three propositions, theoretically, the generator distribution $P_G$ will finally converge to real data distribution $P_{data}$ if $G$ and $D$ have enough capacity, and at each step of training, the discriminator reaches its optimum of the given $G$, and $P_G$ is updated based on the training objective given in eq. (1). The generator of a learned GAN can generate data samples with the same distribution as original data. In other words, the probability distribution learned by a GAN is a good estimation for empirical probability distribution $P_0$. Therefore, solution quality of GAN based data-driven distributionally robust chance constrained programming is guaranteed.

*GAN based distributionally robust chance constrained programming using sample average approximation*

In general, chance constrained programming problem is considered intractable because it is in general impossible to check the feasibility of a given candidate solution, and the feasible region defined by chance constraints is usually nonconvex.[56,57] Therefore, we seek to find approximate solutions to chance constrained programming problems. One common strategy is to find conservative approximations which can be solved efficiently. In other words, the approximate problem could be solved efficiently and can yield feasible solutions or at least feasible with a high probability. SAA is an intuitive and popular method for solving chance constrained programming problems in a data-driven perspective.[37] The basic idea for SAA is to approximate the true probability distribution of uncertain parameters with Monte Carlo sampling based on an empirical probability distribution. The SAA problem (P2) is formulated by replacing the latent uncertain parameters with their empirical counterparts constructed using independent and identically distributed (i.i.d.) data samples of uncertain parameters shown below.

(**P2**)   $\min f(x)$ (21)

$$s.t. \ \frac{1}{N}\sum_{j=1}^{N}\mathbf{1}_{(0,+\infty)}\left(C\left(x,\zeta^{j}\right)\right) \leq \gamma \tag{22}$$

$$x \in X \tag{23}$$

In problem (P2), $N$ is the sample size with each data sample denoted by index $j$, and the indicator function **1** represents whether equation $C$ is great than 0, or in other words whether the constraint is violated or not. $\zeta^j$ are i.i.d. samples of random parameter $\zeta$ generated from empirical distribution $P_0$, and $\gamma$ is the revised risk level, which is usually greater than $\alpha'$. From the previous subsection, the generator distribution $P_G$ for GAN will finally converge to real data distribution and can serve as a good choice of empirical distribution. From existing studies, we know that the generator of GAN can be regarded as a transformation which transforms a random noise variable to a sample under generator distribution.[44,49,50] Thus, from probability theory we can conclude that if the random noise variables which serve as the input of generator network are i.i.d., then the samples generated by the generator network of GAN are i.i.d. Therefore, GAN can be adopted to estimate the empirical probability distribution and generate data samples for SAA simultaneously. If the constraints and objective function of original problem (P0) are linear, the SAA problem (P2) is a mixed-integer linear programming (MILP) problem which can be solved efficiently by branch-and-cut algorithms implemented in off-the-shelf optimization solvers like CPLEX and GUROBI.

From the SAA problem (P2), we can get a candidate solution indicated by $x^*$. The next step is to check its quality or whether it can serve as a good solution to the chance constrained programming problem (P1). Specifically, two questions should be addressed: (1) whether $x^*$ is a feasible solution to the original problem, and (2) how large the optimality gap is.[38] However, it is challenging to get the exact answers to these two questions, so we seek to verify the probability that $x^*$ is a feasible solution to the original problem and estimate the optimality gap.

$$q(x) = \mathbb{P}_{\zeta \sim P}\left(C(x,\zeta) > 0\right) \tag{24}$$

$$\hat{q}_N(x) = \frac{1}{N}\sum_{j=1}^{N}\mathbf{1}_{(0,+\infty)}\left(C\left(x,\zeta^{j}\right)\right) \tag{25}$$

$q(x)$ represents the probability that constraint $C$ is violated and $\hat{q}_N(x)$ is an estimation of $q(x)$. To estimate the probability that constraint $C$ is violated for solution $x^*$, we generate another set of i.i.d. samples from GAN denoted as $\zeta^1, \ldots, \zeta^{N'}$ and estimate $q(x^*)$ by $\hat{q}_{N'}(x^*)$.

$$\hat{q}_{N'}(x^*) = \frac{1}{N'} \sum_{j'=1}^{N} \mathbf{1}_{(0,+\infty)}\left(C(x,\zeta^{j'})\right) \tag{26}$$

$$U_{\beta,N'}(x^*) = \hat{q}_{N'}(x^*) + z_\beta \sqrt{\hat{q}_{N'}(x^*)(1-\hat{q}_{N'}(x^*))/N'} \tag{27}$$

It has been proved that if $U_{\beta,N'}(x^*) \leq \alpha'$, $x^*$ is a feasible solution for the original chance constrained programming problem and $f(x^*)$ is an upper bound at least with probability $1 - \beta$.[38,52,55] $z_\beta$ is a quantile value of standard normal distribution corresponding to $\beta$ with $z_\beta = \Phi^{-1}(1-\beta)$.

Since we have got an upper bound of the optimal value, the next step is to find a lower bound. Here we choose two positive integers $M$ and $N$. $L$ is the largest integer satisfying constraints (28) and (29). $B$ denotes binomial probability distribution quantile value.

$$\theta_N = B(\lfloor \gamma N \rfloor; \alpha', N) \tag{28}$$

$$B(L-1; \theta_N, M) \leq \beta \tag{29}$$

$M$ sets of samples are generated from GAN, and each set of samples has size $N$. By solving the SAA problem (P2), $M$ solutions can be calculated denoted by $f(x_1), \ldots, f(x_M)$. Next, they are sorted in an increasing order shown as eq. (30):

$$f(x_{(1)}) \leq \ldots \leq f(x_{(M)}) \tag{30}$$

It can be proved that the $L$'th smallest value of $f(x)$ is a lower bound with probability at least $1 - \beta$.[38,56]

To summarize, the solution algorithm of GAN based data-driven distributionally robust chance constrained programming is provided in Figure 2. $H$ describes how many times sample average approximation is repeated.

| | GAN based data-driven distributionally robust chance constrained programming solution algorithm |
|---|---|
| 1 | Reformulate problem (P0) to problem (P1) |
| 2 | Train GAN with historical data; |
| 3 | Derive SAA problem (P2) |
| 4 | **while** $gap \geq tol$ **do** |
| 5 |   **for** $h = 1, 2, \ldots, H$, **repeat** |
| 6 |     **for** $m = 1, 2, \ldots, M$, **repeat** |
| 7 |       Generate $\zeta^1, \ldots, \zeta^N$ $N$ data samples by GAN; |
| 8 |       Solve the associated SAA problem (P2) with $N$ scenarios. Denote the solution as $x_m$ and optimal value as $v_m$; |
| 9 |       Generate $\zeta^1, \ldots, \zeta^{N'}$ $N'$ new data samples by GAN; |
| 10 |       Calculate $\hat{q}_{N'}(x_m)$, $U_{\beta,N'}(x_m)$ ; |
| 11 |       If $U_{\beta,N'}(x_m) \leq \alpha'$ , go to next step; If not, skip and go to next iteration; |
| 12 |       $UB = \min\{UB, f(x_m)\}$ denoted as $g_h$; |
| 13 |     **End** |
| 14 |     Sort the $M$ optimal values and find the $L$th smallest one as lower bound $LB$ denoted as $vl_h$; |
| 15 |   **End** |
| 16 |   Calculate $v = H^{-1}\sum_h vl_h$ as the lower bound; |
| 17 |   Calculate $g = \min_h g_h$ as the upper bound; |
| 18 |   $gap = (g-v)/g$ |
| 19 | **End** |
| 20 | **return;** |

Figure 2 Solution algorithm for GAN based data-driven distributionally robust chance constrained programming.

*Summary*

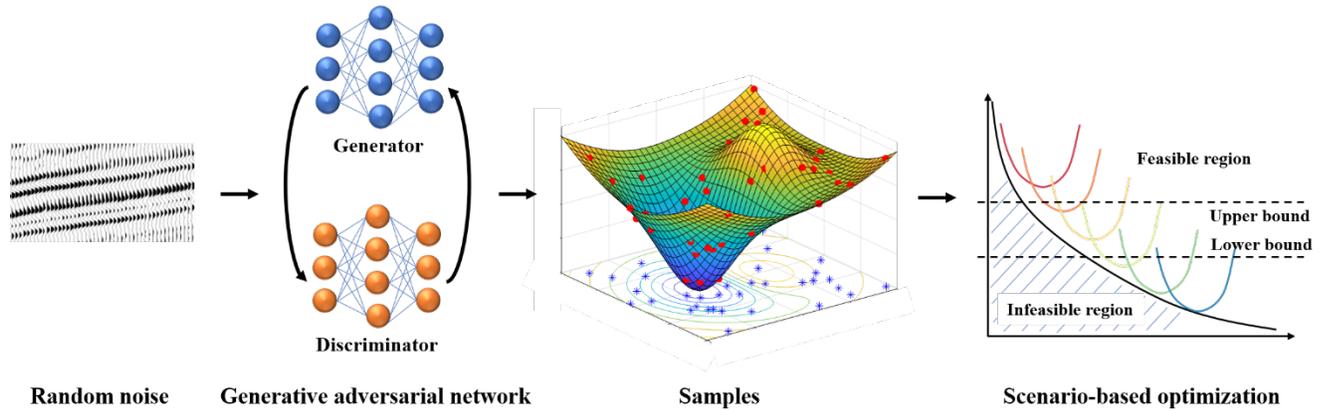

Figure 3 GAN based data-driven distributionally robust chance constrained programming framework.

In this section, we introduce the whole framework of GAN based data-driven distributionally robust chance constrained programming. The proposed framework can be illustrated through Figure 3. The reformulation steps and SAA method do not require linear objectives and constraints and thus the whole framework still works for nonlinear problems. Compared with other probability distribution estimation methods, GAN-based framework has three advantages.

- First, conventional density estimation methods usually adopt approximation methods and prior assumptions, which can make the optimal solutions to the data-driven optimization problems sensitive and highly dependent to the approximation methods and assumptions. In the contrary, GAN does not require a priori approximation or assumption, and can be trained in an unsupervised way so that the dependency and sensitivity to the assumptions can be avoided for optimal solutions to data-driven optimization problems.[50]
- Second, since the generator and discriminator of a GAN are usually deep neural networks, more complicated probability distributions can be learned, and even some hidden modes can be found compared with conventional probability distribution estimation methods. GAN can efficiently and accurately model uncertainty,[58] and thus increasing the solution quality to data-driven optimization problems.
- Third, conventional density estimation methods cannot model uncertainty accurately and correctly if there are missing data in the data set. GAN can make compensation to the missing information in the available data set, which can reduce the influence of missing data in uncertainty modelling.[59] In other words, GAN can utilize the remaining data in the data set to estimate the missing data distribution. In this way, GAN can lead to more accurate optimal solutions to data-driven optimization problems. This is an important advantage and will be illustrated in the next section.

*Illustrative example*

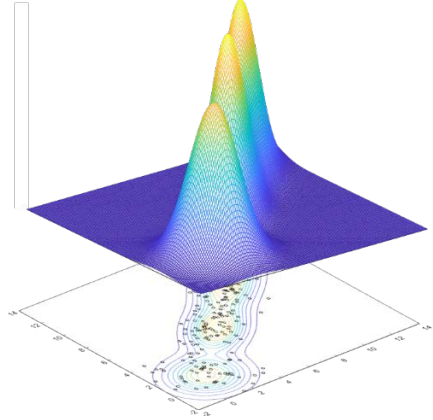

Figure. 4 True probability distribution for generating data samples with density contour.

To further illustrate the advantage of the proposed framework, we present an illustrative example. In this example, we compare the distribution approximation results of GAN with results from anther non-parametric distribution estimation method, kernel density estimation (KDE). The results can clearly show that GAN can better approximate the unknown distribution especially in the case when there is missing data in the original data set. Suppose the data set has a true probability distribution of a Gaussian mixture distribution as shown in Figure 4 and samples are generated according to this true distribution. The whole initial data set is shown in Figure 5 (a). However, there exists missing data in the available data set due to various reasons. In this case, data in the rectangular area $\{(x, y): 6 \leq x \leq 8, 6 \leq y \leq 8\}$ is missing with a ratio of 21/200, as shown in the red rectangular area in Figure 5 (a). Therefore, we can only utilize the remaining data to estimate the empirical probability distribution and then generate data samples. For comparison, KDE and GAN are both employed to model the probability distribution with the remaining data. Figure 5 (b) shows the density contour of the probability distribution estimated by KDE. It can be found that in the red rectangular area $\{(x, y): 6 \leq x \leq 8, 6 \leq y \leq 8\}$, the probability distribution estimated by KDE results in a relatively low density compared with the density of the original probability distribution. Figures 5 (c) and (d) show the data samples generated by KDE and GAN, respectively. It can be found that KDE generates data samples in the rectangular area $\{(x, y): 6 \leq x \leq 8, 6 \leq y \leq 8\}$ with a low probability, as shown by the ratio of points in the red rectangular area of 3/200. In contrast, GAN generates data samples in the rectangular area with a relatively higher probability, with the ratio of points in the rectangular area of 20/200. Thus, KDE does not model uncertainty as accurately as GAN, if the original data set has missing data. However, to some extent, GAN can make compensation to the missing information in the original data set. Besides, we have tried uniform distribution, one of the traditional approximate distributions, to estimate the unknown distribution and generate data samples. However, uniform distribution leads to even worse results than the results from KDE. Therefore, we do not include the results from uniform distribution in the paper.

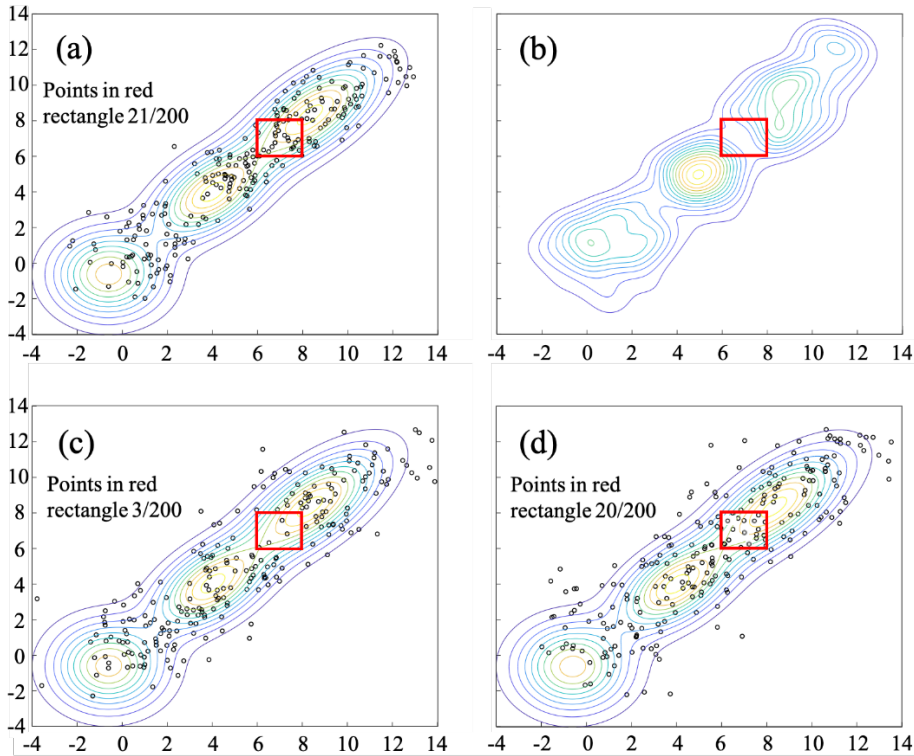

Figure. 5 Illustrate example: (a) original data set with points in red rectangular area, 21/200; (b) density estimated by KDE with remaining data; (c) data samples generated by KDE with points in red rectangular area, 3/200; (d) data samples generated by GAN with points in red rectangular area, 20/200.

## Motivating Example

In this section, we present a small numerical example to illustrate the proposed GAN based data-driven distributionally robust chance constrained programming framework. In this example, the convergence of GAN with different initial points is demonstrated. The example problem (P3) is given as follows:

(**P3**) $\quad \min x_1 + x_2 + x_3$ $\hfill (31)$

$$s.t. \inf_{P \in D} \mathbb{P}_{\xi \sim P} \begin{pmatrix} x_1 + x_2 + \xi_1 \leq 10 \\ x_2 + x_3 + \xi_2 \leq 11 \\ x_1 + x_3 + \xi_3 \leq 12 \end{pmatrix} \geq 1 - \alpha = 1 - 10\% \quad (32)$$

$$D = \{P \in M : \phi(f, f_0) \leq d = 10\%, f = dP/d\xi\} \quad (33)$$

$$x_1, x_2, x_3 \geq 0 \quad (34)$$

where $x_1$, $x_2$, $x_3$ are decision variables and $\xi_1$, $\xi_2$, $\xi_3$ are uncertain parameters with joint probability distribution $P$ in the ambiguity set $D$. To solve this problem, we first define $g(x)$ as in eq. (10) and calculate the modified risk level $\alpha$' as shown by eq. (35).

$$\alpha' = \alpha - \frac{d}{2} = 10\% - \frac{10\%}{2} = 5\% \quad (35)$$

The distributionally robust chance constrained programming problem (P3) is reformulated into a conventional chance constrained programming problem shown as problem (P4),

(**P4**) $\quad \min x_1 + x_2 + x_3$ $\hfill (36)$

$$s.t. \mathbb{P}_{\xi \sim P_0} \begin{pmatrix} x_1 + x_2 + \xi_1 \leq 10 \\ x_2 + x_3 + \xi_2 \leq 11 \\ x_1 + x_3 + \xi_3 \leq 12 \end{pmatrix} \geq 95\% = 1 - 5\% \quad (37)$$

$$x_1, x_2, x_3 \geq 0 \quad (38)$$

where $P_0$ is the empirical probability distribution of uncertain parameters $\xi_1$, $\xi_2$, $\xi_3$. Next, corresponding SAA problem (P5) is derived as shown below:

(**P5**) $\quad \min x_1 + x_2 + x_3$ $\hfill (39)$

$$s.t.\ x_1 + x_2 - A_j \cdot \xi_1^j \leq 10 - \xi_1^j,\ \forall j \tag{40}$$

$$x_2 + x_3 - A_j \cdot \xi_2^j \leq 11 - \xi_2^j,\ \forall j \tag{41}$$

$$x_1 + x_3 - A_j \cdot \xi_3^j \leq 12 - \xi_3^j,\ \forall j \tag{42}$$

$$\frac{1}{N}\sum_{j=1}^{N} A_j \leq 5\% \tag{43}$$

$$x_1, x_2, x_3 \geq 0, A_j \in \{0,1\} \tag{44}$$

The empirical distribution of uncertain parameters $\xi_1$, $\xi_2$, $\xi_3$ is estimated by GAN. Figure 6 illustrates the process of GAN training. Two random noise distributions regarded as inputs of GAN are presented as initial guess distributions (in blue color) in Figures 6 (a) and (b). The target distributions for Figures 6 (a) and (b) are the same. As stated previously, with sufficiently many training epochs, the generator distributions corresponding to Figures 6 (a) and (b) will finally converge to the same target distribution. Intuitively, the training process can be stopped, if the two generator distributions in Figures 6 (a) and (b) are close enough. In this motivating example, we adopt the center-to-center distance to evaluate the similarity of two distributions for stopping criteria. If the distance between the centers of two point distributions generated by these two generator distributions corresponding to Figures 6 (a) and (b) is less than or equal to 5% of the distance between the centers of two initial guess distributions, the GAN training process will be stopped. The center-to-center distance of the two initial guess distributions is 4.33. After training with 50 epochs, GAN generates data samples as the epoch 50 distributions (in green color) in Figures 6 (a) and (b) with a center-to-center distance of 0.97 (22.4% of the center-to-center distance between the two initial guess distributions). After 100 epochs, GAN generates data samples as the epoch 100 distributions (in red color) in Figures 6 (a) and (b), which are highly consistent with the target distribution, with a center-to-center distance of 0.21 (4.8% of the center-to-center distance for the two initial guess distributions). With the increase of training epochs, the point groups generated by the generator distributions in Figures 6 (a) and (b) gradually overlap. Specifically, the epoch 50 distributions in Figures 6 (a) and (b) only have some small overlap, while the epoch 100 distributions in Figures 6 (a) and (b) overlap significantly. With 100 training epochs, the center-to-center distance between the resulting distributions is less than 5% of the distance between the centers of two initial guess distributions. Therefore, we can stop training at 100 epochs.

In addition, we can solve the corresponding SAA problems with these three sample distributions in Figures 6 (a) and (b) and get the optimal objective function values. From the previous section, when $\beta = 5\%$, eq. (31) and (32) can be satisfied with $N = 200$ and $M = 20$. The optimality tolerance calculated by $(g - v)/g$ shown in Figure 2 is set to be 1%. To determine the value of $H$ in the solution algorithm shown in Figure 2, we first set $H = 10$ and the resulting optimality gap is greater than 1%, which indicates $H = 10$ cannot satisfy the optimality tolerance requirement. Next, $H$ is increased to 20 and the optimality tolerance can be satisfied. Thus, $H$ is set to 20 for this motivating example. According to the initial guess distributions in Figures 6 (a) and (b), the optimal objective functions calculated from the solution algorithm are 12.1 and 9.5, respectively. Similarly, the optimal objective functions corresponding to the epoch 50 distributions and the epoch 100 distributions can be obtained. The optimal objective functions corresponding to epoch 50 and 100 distributions in Figure 6 (a) are 13.5 and 14.2, respectively. The optimal objective functions corresponding to epoch 50 and 100 distributions in Figure 6 (b) are 12.3 and 14.5, respectively. Ideally, the SAA problems corresponding to Figures 6 (a) and (b) will result in almost the same optimal solution with adequate data samples and same data distribution. Although the epoch 100 distributions in Figures 6 (a) and (b) are close enough, there exist small differences which could cause differences in data samples. Moreover, even with the same sample distribution, the number of data samples should be limited when solving an SAA problem due to computational tractability. There may be minor difference between the optimal objective functions corresponding to Figures 6 (a) and (b) with epoch 100 distributions. With different initial points, the GAN will converge to the same distribution.

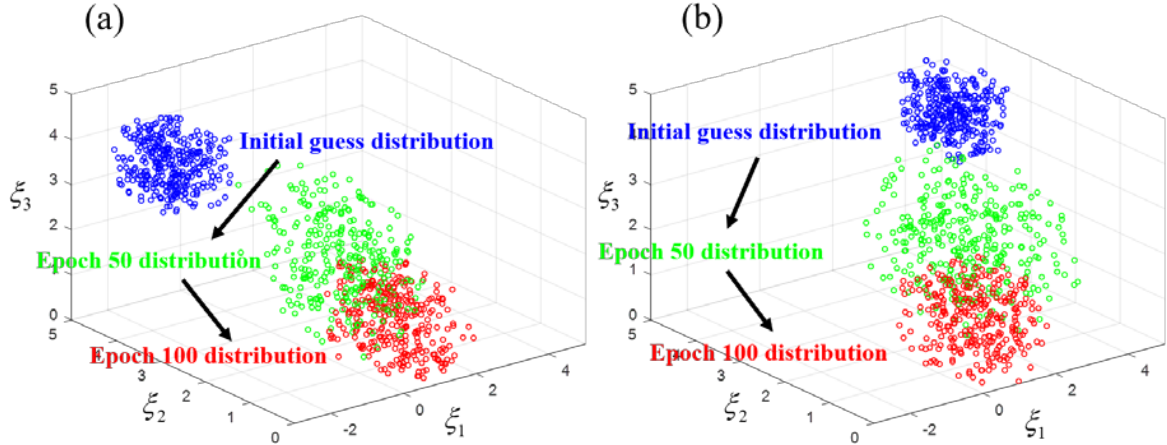

Figure 6  Illustration of GAN training process with initial guess distributions (in blue color), epoch 50 distributions (in green color) and epoch 100 distributions (in red color) from two approaches. Cases (a) and (b) have different initial guess distributions as input to GAN and the same target distribution. Training stops at 100 epochs when the distance between the centers of two generator distributions (0.21) is less than or equal to 5% of the distance between the centers of two initial guess distributions (4.33).

## Application to Supply Chain Optimization under Uncertainty

Supply chain optimization is an important area of practical importance.[60-62] Many existing studies adopt data-driven optimization frameworks to optimize the facility locations, inventory levels, network structures and transportations in supply chains.[63,64] In practice, companies are faced with massive influx of big data with complex correlations from different participants in a supply chain.[65] Since uncertainty plays a key role in the data-driven optimization frameworks, it is important to efficiently and correctly model uncertainty with the complex correlated data in supply chain optimization problems.[39] Thus, deep learning based data-driven optimization methods fit well with supply chain optimization problems.

In this section, we apply the proposed data-driven distributionally robust chance constrained programming framework to supply chain optimization under uncertainty. The problem statement, distributionally robust chance constrained programming model formulation, solution algorithm and two-stage stochastic programming model formulation are given as follows. To illustrate the applicability, a case study of a spatially explicit biofuel supply chain in Illinois is presented.

*Problem statement*

In this work, we address the optimal design and operations of supply chains under demand uncertainty. The supply chain can be represented as a three-echelon network, including manufacturing facilities, customers and raw material suppliers. We are given a set of candidate sites for manufacturing facility location, a set of customers and a set of raw material suppliers. A finite planning horizon is considered, and it is partitioned into a set of time periods. The duration of each time period is known.

For each manufacturing facility, the capital cost is considered as a piecewise linear function of production capacity. A set of production capacity levels are given with the upper and lower bounds of each production capacity level. The unit capital cost for facility location and the conversion rate from the raw material to the product are also given.

For each customer, the demand in each time period is treated as an uncertain parameter with an unknown ambiguous probability distribution. The unknown probability distribution is considered belonging to a set of possible probability distributions estimated using existing historical data.[39,66] For each raw material supplier, the unit raw material procurement cost and the maximum supply of raw materials are considered as given.

The transportation flows include those from suppliers to manufacturing facilities and from manufacturing facilities to customers. The unit transportation cost for product between each candidate facility location and each customer, as well as the unit transportation cost for raw material between each candidate facility location and each supplier, are given.[67]

The inventory level in each manufacturing facility of current time period is determined by transportation flows, production level and inventory level of the previous time period. The inventory degradation after each time period is a fixed percentage of the total inventory level at the previous time period and the inventory degradation rate is considered as given.[68]

In the supply chain optimization problem, facility locations, production capacities, external supplies, transportation, procurement and inventory decisions need to be determined to satisfy the demand requirements with at least a relatively high probability in the worst case. The worst case is taken over a set of possible probability distributions of customer demand. In other words, the demand requirements should be satisfied for a set of possible probability distributions with at least a pre-defined probability. The objective for the problem is to minimize the total cost, including capital cost, operating cost, transportation cost, procurement cost, inventory cost and external product cost.

*Distributionally robust chance constrained programming model formulation*

According to the previous problem statement, we can propose the distributionally robust chance constrained programming model for supply chain optimization. In the model formulation, there are $I$ candidate sites for manufacturing facility location (indexed by $i = 1, \ldots, I$), $K$ suppliers for raw material (indexed by $k = 1, \ldots, K$), $J$ customers (indexed by $j = 1, \ldots, J$), $R$ production capacity levels for a manufacturing facility (indexed by $r = 1, \ldots, R$) and $T$ time periods (indexed by $t = 1, \ldots, T$). The problem is formulated in discrete time with finite time horizon.

$$\sum_r Y_{ir} \leq 1, \ \forall i \tag{45}$$

$$CF_i = \sum_r \left( cm_{r-1} \cdot Y_{ir} + \frac{(CA_{ir} - pm_{r-1} \cdot Y_{ir})(cm_r - cm_{r-1})}{pm_r - pm_{r-1}} \right), \ \forall i \tag{46}$$

$$pm_{r-1} \cdot Y_{ir} \leq CA_{ir} \leq pm_r \cdot Y_{ir}, \ \forall i, r \tag{47}$$

$$CO_i = vv_i \cdot \sum_r CA_{ir}, \ \forall i \tag{48}$$

The facility location constraints are given as constraints (45) – (48). Binary variable $Y_{ir}$ denotes the facility location and production capacity level decision. $Y_{ir}=1$ represents the manufacturing facility with capacity level $r$ is selected to build at candidate site $i$. Constraint (45) enforces that at most one capacity level can be selected in each candidate site. $CA_{ir}$ denotes the production capacity in each time period for manufacturing facility $i$ in capacity level $r$. The capital cost of facility location is a piecewise linear function with respect to capacity $CA_{ir}$. Constraint (46) calculates the capital cost for facility location at candidate site $i$, denoted by $CF_i$. Constraint (47) sets the lower and upper bounds for capacity at each capacity level. $pm_r$ denotes the upper bound of capacity at capacity level $r$, and $cm_r$ indicates the reference capital cost corresponding to reference capacity $pm_r$. Constraint (48) calculates the operating cost at candidate site $i$, denoted by $CO_i$. Parameter $vv_i$ represents the unit production cost in manufacturing facility $i$.

$$\sum_j X_{ijt} \leq \sum_r CA_{ir}, \ \forall i, t \tag{49}$$

$$\sum_i Z_{kit} \leq bm_{kt}, \ \forall k, t \tag{50}$$

The transportation constraints are given as constraints (49) – (50). Constraint (49) enforces that the operating level of each facility should not exceed the established capacity. $X_{ijt}$ denotes the transportation amount of product from facility $i$ to customer $j$ at time period $t$. Constraint (50) describes the supply side of the supply chain. $Z_{kit}$ is the transportation amount of raw materials from supplier $k$ to facility $i$ at time period $t$, and $bm_{kt}$ is the maximum raw material supply in supplier $k$ at time period $t$.

$$(1-\eta) \cdot \left( I_{it} + \sum_k Z_{kit} - \sum_j X_{ijt} / \beta \right) = I_{it+1}, \forall i, 1 \leq t < T \tag{51}$$

$$(1-\eta) \cdot \left( I_{iT} + \sum_k Z_{kiT} - \sum_j X_{ijT} / \beta \right) = I_{i1}, \forall i \tag{52}$$

The inventory constraints are given as constraints (51) – (52). Constraint (51) calculates the inventory level at each facility in each time period. $I_{it}$ denotes the inventory level at facility $i$ in time period $t$. $\beta$ is the conversion rate from raw material to product and $\eta$ is the raw material inventory degradation rate for each time period. Constraint (52) describes the inventory levels in the beginning of the first time period equals to that at the end of the last time period.

$$\inf_{P \in D} \mathbb{P}_{d_{jt} \sim P} \left( \sum_i X_{ijt} + L_{jt} \geq d_{jt}, \forall j, t \right) \geq 1 - \alpha \tag{53}$$

$$D = \{ P \in M : \phi(f, f_0) \leq d \} \tag{54}$$

The customer demand distributionally robust chance constraints are given in constraints (53) – (54). The customer demand values at each time period are considered as uncertain parameters denoted by $d_{jt}$ with an unknown joint probability distribution $P$ belonging to the ambiguity set $D$. $L_{jt}$ denotes the external supply of product at customer $j$ in time period $t$. From another perspective, the value of external supplies can quantify the feasibility or infeasibility of the solution with the case without external supplies. Larger external supplies imply that the corresponding result will have a higher chance to be infeasible. Constraint (53) describes the demand side of the supply chain and is the distributionally robust chance constraint. Constraint (53) states that the summation of delivered product and external product at each customer in each time period should be greater than or equal to the corresponding demand by at least a pre-specified probability $1 - \alpha$ in the worst case. The worst case is taken over ambiguity set $D$. Constraint (54) defines the ambiguity set through $\phi$-divergence and the range of the ambiguity set is denoted by $d$. Two approaches can be adopted to consider more "worst case" scenarios or "worst case" scenarios with higher risks, increasing the range of the ambiguity set $d$ and increasing the pre-defined probability of the chance constraints $1 - \alpha$.

$$\min \sum_i (CF_i + CO_i) + pi \sum_{i,t} I_{it} + \sum_{i,j,t} c_{ij} X_{ijt} + \sum_{k,i,t} a_{ki} Z_{kit} + pex \sum_{jt} L_{jt} \tag{55}$$

The objective function is described by eq. (55) that calculates the total cost including capital cost, operating cost, transportation cost, procurement cost, inventory cost and external product cost. *pi* is the unit inventory cost and *pex* is the unit cost for external product. $c_{ij}$ is the unit transportation cost for product from facility $i$ to customer $j$. $a_{ki}$ is the unit transportation cost for raw material from supplier $k$ to facility $i$.

For summary, the proposed distributionally robust chance constrained programming model for supply chain optimization can be described as the following problem (P6).

**(P6)** min objective function given in eq. (55)
s.t. facility location constraints (45) – (48)
transportation constraints (49) – (50)
inventory constraints (51) – (52)
customer demand distributionally robust chance constraints (53) – (54)

*Solution algorithm*

To solve the proposed distributionally robust chance constrained programming problem (P6), we follow the steps described in the GAN-based Data-Driven Distributionally Robust Chance Constrained Programming Section. First, to construct the ambiguity set, $g(x)$ is set as follows to calculate $\phi$-divergence for computational efficiency:[35]

$$g(x) = |x-1| \tag{56}$$

Second, we train the GAN with historical demand data and reformulate the distributionally robust chance constrained programming problem (P6) into the corresponding conventional chance constrained programming problem (P7) with a modified risk level $\alpha'$:

$$\alpha' = \alpha - \frac{d}{2} \tag{57}$$

Problem (P7) is obtained by replacing the distributionally robust chance constraints (53) – (54) with the corresponding conventional chance constraint (58) under empirical probability distribution $P_0$ and modified risk level $\alpha'$.

**(P7)** min objective function given in eq. (55)
*s.t.* constraints (45) – (52)

$$\mathbb{P}_{d_{jt} \sim P_0}\left(\sum_i X_{ijt} \geq d_{jt},\ \forall j,t\right) \geq 1-\alpha' \tag{58}$$

Next, the conventional chance constrained programming problem (P7) is converted to the SAA problem (P8) shown as follows:

(**P8**) min objective function given in eq. (55)

  *s.t.* constraints (45)– (52)

$$\sum_i X_{ijt} + d_{jtn} \cdot A_n \geq d_{jtn},\ \forall j,t,n \tag{59}$$

$$\sum_n A_n \leq N \cdot \gamma \tag{60}$$

In problem (P8), constraints (59) – (60) are the sample average approximation to the chance constraint (58) in problem (P7). $d_{jtn}$ are i.i.d. samples of uncertain parameter $d_{jt}$. Binary 0-1 variable $A_n$ indicates if the demand constraint is violated. $N$ is the number of data samples with each data sample indexed by $n$. $\gamma$ is the revised risk level, and $\gamma$ is usually greater than $\alpha'$ with the following relationship:

$$\gamma > \alpha' = \alpha - \frac{1}{2}d \tag{61}$$

$\hat{q}_N(x)$ for solution validation is given in eq. (62).

$$\hat{q}_{N'}(x^*) = \frac{1}{N'} \sum_{n'=1}^{N} \mathbf{1}_{\left(\sum_i X_{ijt} < d_{jm'},\ \forall j,t\right)} \tag{62}$$

where **1** stands for the indicator function. Finally, following the algorithm given in Figure 2, we can get the optimal solution of the GAN based distributionally robust chance constrained programming problem.

## *Case study of spatially explicit biofuel supply chain in Illinois*

  The applicability of the proposed approach is illustrated through a county-level case study of a spatially explicit biofuel supply chain in Illinois. Mathematical programming methods have been extensively employed for biofuel and bioenergy supply chain optimization.[69,70] A typical biofuel supply chain contains a three-echelon superstructure including biomass suppliers, biorefineries and biofuel customers. This biofuel supply chain optimization problem needs to determine the optimal biorefinery locations, production capacities, biomass transportation from biomass suppliers to biorefineries, biofuel transportation from biorefineries to biofuel customers and inventory levels.[71] The biofuel product is fuel ethanol, and the feedstock is corn stover. We consider 25 counties in Illinois with largest population and biofuel consumption as potential biofuel customers, and 15 counties with moderate corn stover yield and biofuel consumption as candidate sites for biorefinery locations. Biofuel from other states is treated as external biofuel supply. Biomass can be acquired from counties within Illinois and from other states. The top 25 counties in Illinois with highest corn stover yields are considered as potential biomass suppliers; there is also a "dummy" biomass supplier with much higher transportation cost to represent the biomass suppliers from other states. Therefore, totally 26 biomass suppliers are taken into account in our problem. To account for the seasonality in biomass supply, we consider a time horizon of one year and each time period represents one month.[68,72] Based on existing studies, corn stover can only be harvested in October and November.[72] Model parameters and techno-economic data are obtained from the literature and various public sources.[72-76]

We first use the customer demand data from 1983 to 2018 as inputs to train the GAN to generate data samples for the SAA method.[72] GAN is implemented and trained with PyTorch 1.1.[77] In the optimization model, the risk level $\alpha$ is set to be 10% and the range of ambiguity set $d$ equals to 0.1. Thus, the modified risk level $\alpha$' equals to 5%. The average tortuosity factor for biomass and biofuel transportation is taken as 1.6.[78] The corresponding SAA MILP problem (P8) is solved with Gurobi 8.1. All instances are solved using a PC with an Intel Core i7-6700 CPU at 3.40 GHz and 32.00 GB RAM. Optimization models and solution procedures are coded in Julia v1.1 with JuMP package.[79] The optimality tolerance for Gurobi 8.1 solver is set to 0.1%. The optimality tolerance for solution algorithm is set to 1%. With 50 samples, the corresponding SAA MILP problem (P8) has 95 binary variables, 9,735 continuous variables and 15,808 constraints. The whole solution algorithm takes around 80,283 CPUs (around 22.3 hours).

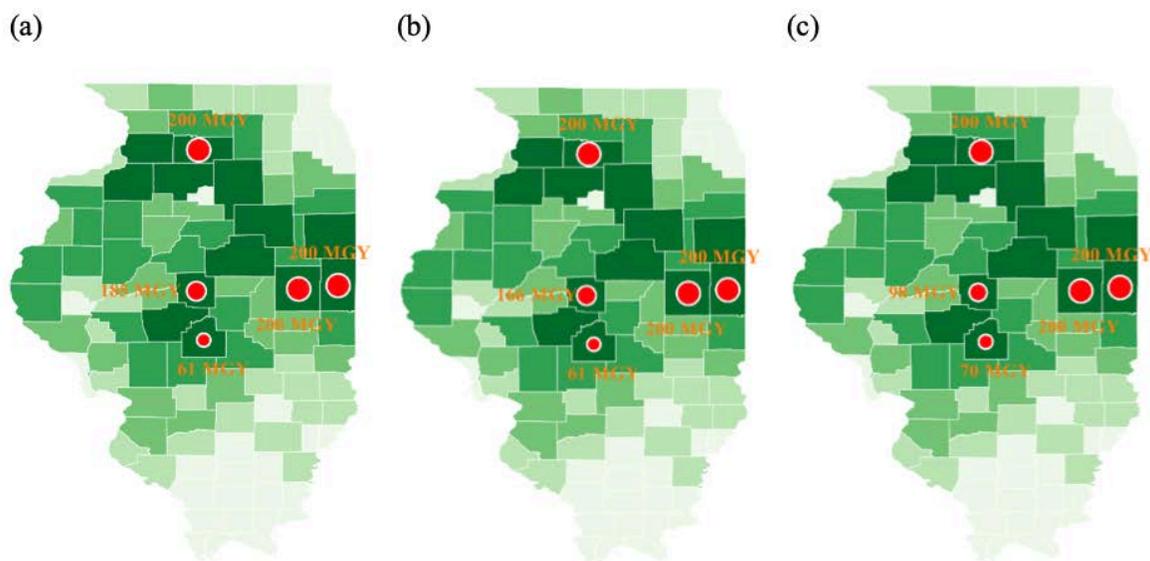

Figure 7 Optimal designs of the biofuel supply chain from different optimization models: (a) distributionally robust chance constrained programming; (b) two-stage stochastic programming; (c) deterministic.

The optimal design of the biofuel supply chain network from distributionally robust chance constrained programming model is shown in Figure 7 (a). Note that the background of the map indicates corn stover yield for each county in Illinois. The optimal total cost is $4,228 MM. The biorefineries in Champaign county, Lee county and Vermilion county are all at the maximum capacity of 200 MM gallons/year. The biorefineries in Christian county and Logan county have smaller capacities of 61 MM gallons/year and 185 MM gallons/year, respectively. From the optimal supply chain design, it can be found that the biorefineries tend to be built in counties with high corn stover yields. The reason is that compared with biofuel transportation cost, biomass acquisition and transportation cost makes up a relatively large part in total cost, which will be clearly illustrated in the cost breakdown.

To take a further look at the optimal solution, we show the cost breakdown of the optimal solution to distributionally robust chance constrained programming problem with GAN generated samples in Figure 8. Note that the biomass acquisition and biomass transportation cost are merged together in the cost breakdown profile. The largest cost comes from biomass acquisition and biomass transportation accounting for 47% of the total cost. This indicates that the purchase price of biomass and the biomass-related logistics costs are considerable in the total cost. The second largest cost comes from capital investment for building biorefineries, which accounts for 27%. The biomass inventory holding cost accounts for 12% of the total cost. This is because biomass can only be harvested in limited months and a large quantity of biomass must be stored over time to minimize material degradation throughout the year. External biofuel supply accounts for 12% of the total cost. In contrast, the biofuel-related logistics costs are much lower, because ethanol has a higher density and is easier for transportation. Figure 9 shows the transportation flows of biofuel product, where each link represents a biofuel transportation flow from a biorefinery to a customer county.

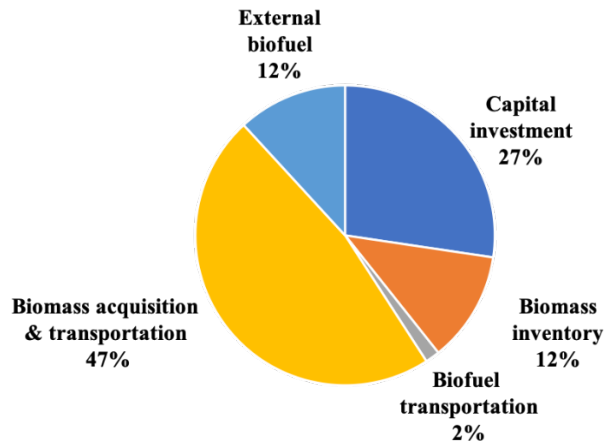

Figure 8 Cost breakdown of the optimal solution to the distributionally robust chance constrained programming problem with GAN generated samples.

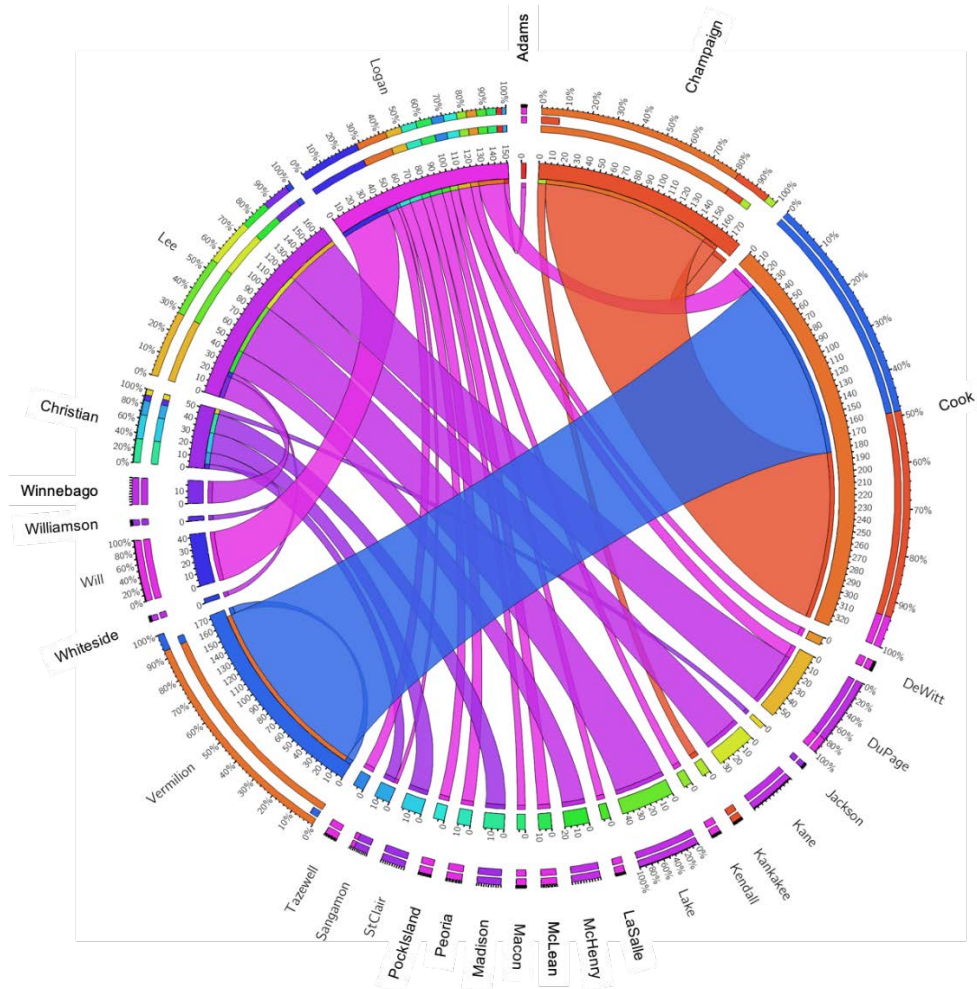

Figure 9 Biofuel product transportation flows of the optimal solution to the distributionally robust chance constrained programming problem with GAN generated samples.

For comparison, the deterministic problem and the two-stage stochastic programming problem are solved. The optimal supply chain designs corresponding to the two-stage stochastic programming problem and deterministic problem are presented in Figures 7 (b) and (c), respectively. The optimal expected total cost of two-stage stochastic programming model is $4,023 MM. The deterministic model results in a total cost of $3,716.

Furthermore, we make a comprehensive comparison between optimal solutions from distributionally robust chance constrained programming model, deterministic model and two-stage stochastic programming model. The biomass inventory profiles and the biofuel production level profiles corresponding to these optimization models are presented in Figure 10. Each column in the inventory profile represents the level of biomass inventory in a certain month. The biomass inventory profiles from these three optimization models all have similar trends. The maximum biomass inventory level is in December, and the inventory level keeps decreasing from December to October in the second year until it reaches zero in October. Because biomass can only be harvested in October and November, large amount of biomass tends to be purchased in November instead of October in order to reduce the inventory cost. However, the production levels from three optimization models have different trends and fluctuate with the customer demand.

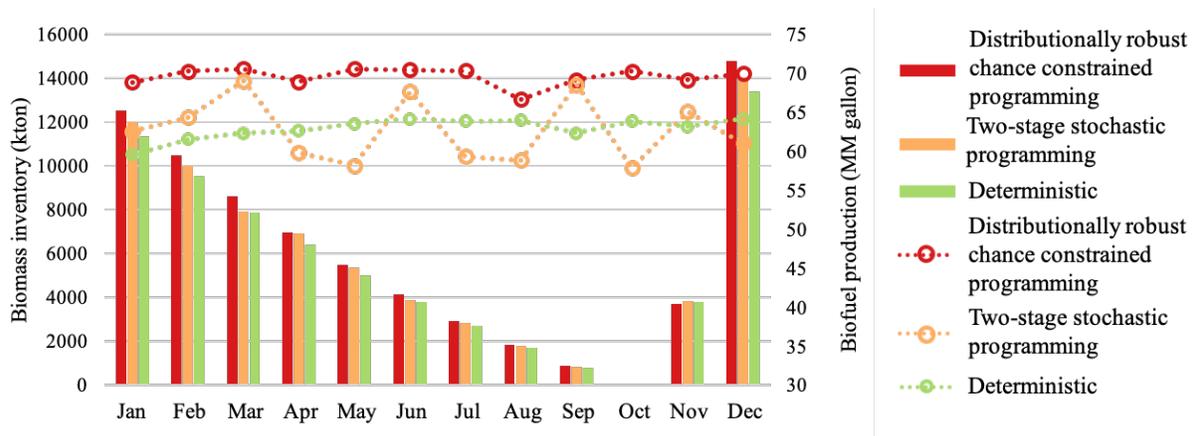

Figure 10 Total biomass inventory and biofuel production profiles in different solutions.

The optimal solution to the deterministic problem has the lowest total cost, production capacity, total biomass inventory level and biofuel production level. However, the corresponding supply chain design could be suboptimal or even infeasible, because it does not consider uncertainty explicitly. The optimal solution to the distributionally robust chance constrained programming problem leads to the highest total cost, production capacity, total biomass inventory level and biofuel production level. The optimal solution to the two-stage stochastic programming problem yields an intermediate total cost, production capacity, total biomass inventory level and biofuel production level. Therefore, the optimal solution to the deterministic problem can be considered as the cost-effective solution and the optimal solution to the distributionally robust chance constrained programming problem can be regarded as the risk-averse solution. In practice, the deterministic solution could be suboptimal or even infeasible due to uncertainty, because it does not account for uncertainty explicitly. Distributionally robust chance constrained programming model accounts for robustness in uncertainty and considers the worst case, which would lead to a higher total cost.

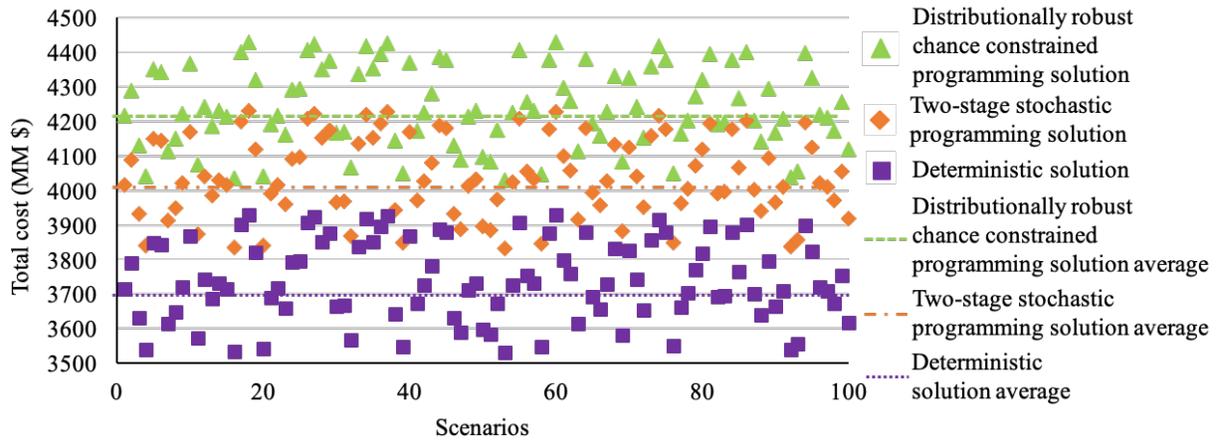

Figure 11 Simulation results for different supply chain designs with data samples from historical realization data.

To further demonstrate the performances of different optimization models for supply chain designs in an uncertain environment, we perform a simulation with 100 data samples from historical realization data. The results are presented in Figure 11, where the X-axis and Y-axis represent the number of scenarios and optimal total cost of the whole supply chain network corresponding to three supply chain designs, respectively. The green triangles denote the total cost values corresponding to distributionally robust chance constrained programming model. These solutions are obtained by solving a deterministic problem with fixed supply chain design as shown in Figure 7 (a). We observe that these solutions indicate a higher average total cost as the green dash line and thus can be considered as risk-averse solutions with lower risk level. The purple squares denote the total cost values corresponding to deterministic models, which are obtained by solving a deterministic problem for each uncertainty realization with the fixed supply chain design as shown in Figure 7 (c). These solutions indicate the lowest average total cost as the purple dash line. The orange diamonds denote the total cost values corresponding to two-stage stochastic programming models, which are obtained by solving a deterministic problem for each uncertainty realization with the fixed supply chain design as shown in Figure 7 (b). These solutions indicate an intermediate average total cost as the orange dash line. Note that for some scenarios with relatively high demand, external biofuel product is needed to satisfy customer demand. Detailed external biofuel product supply is presented as Figure 12 that includes the total number of scenarios using external biofuel and average external biofuel supply cost in the optimal solutions determined by the three approaches. For supply chain design corresponding to distributionally robust chance constrained programming model, external biofuel is needed only under a small number of scenarios with a low average external biofuel supply cost. However, for supply chain design corresponding to the deterministic optimization model, external biofuel is needed with a relatively large number of scenarios causing a high average external biofuel supply cost. Thus, the simulation results explicitly indicate the trade-off between total cost and risk level, where a lower risk level leads to a higher total cost.

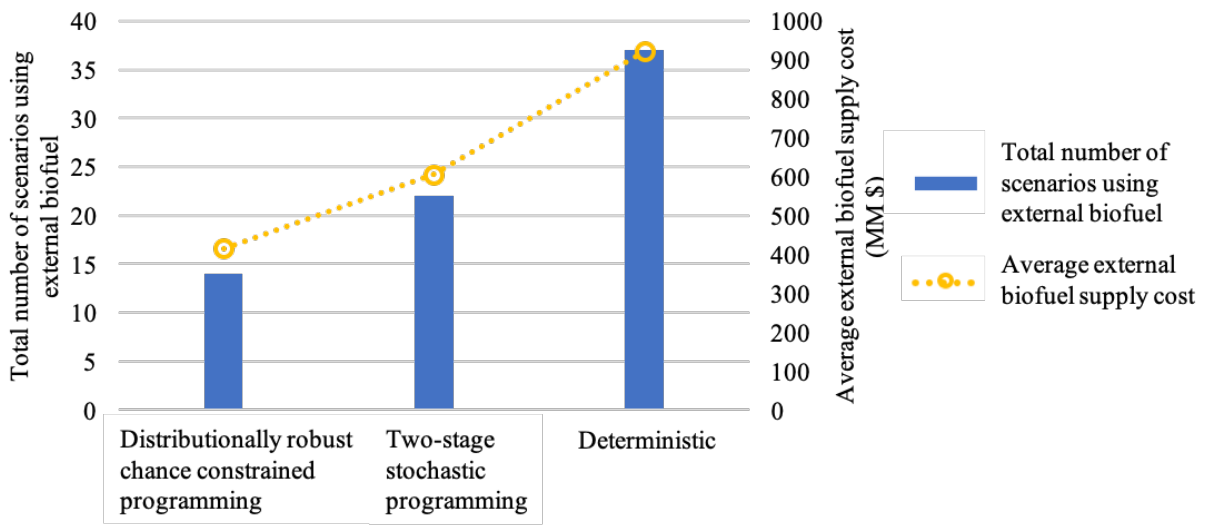

Figure 12 External biofuel product supply including the total number of scenarios using external biofuel and average external biofuel supply cost in the optimal solutions of the three approaches.

## Conclusion

In this work, we proposed a GAN based data-driven distributionally robust chance constrained programming framework. GAN was employed to model uncertainty in data-driven distributionally robust chance constrained programming problems from historical data in a nonparametric way. Distributionally robust chance constrained programming problems were reformulated into conventional chance constrained programming problems, and finally solved by SAA method. The required data samples in SAA method were generated by GAN in an end-to-end way through the differentiable networks. Thus, GAN could extract distributional information from historical data and generate data samples simultaneously. We then applied this framework to a three-echelon supply chain optimization problem under demand uncertainty. The demand for each customer in each time period was treated as an uncertain parameter with an unknown ambiguous probability distribution. The unknown probability distribution was described with an ambiguity set constructed using existing historical data. Facility locations, production capacities, external supplies, transportation, procurement and inventory decisions were determined to satisfy the demand constraints with at least a pre-defined probability in the worst case. The worst case was taken over all the possible probability distributions in the ambiguity set. The objective for the problem was to minimize the total cost including capital cost, operating cost, transportation cost, procurement cost, inventory cost and external product cost. The applicability of the proposed approach was illustrated through a county-level case study of a spatially explicit biofuel supply chain in Illinois. For comparison, we also calculated the optimal solutions to the corresponding deterministic problem and two-stage stochastic programming problem. Results show that the optimal solution of the distributionally robust chance constrained programming problem could be regarded as the risk-averse solution which would lead to a lower risk level and a higher total cost.

## Nomenclature

**Sets/Indices**

| | |
|---|---|
| $D$ | ambiguity set for the distributionally robust chance constrained programming problem |
| $I$ | set of candidate site of facility location indexed by $i$ |
| $J$ | set of customers indexed by $j$ |
| $K$ | set of suppliers indexed by $k$ |
| $M$ | set of all probability distributions |
| $R$ | set of production capacity levels indexed by $r$ |
| $S$ | probability distribution of uncertain parameters in the two-stage stochastic programming problem with scenarios indexed by $s$ |
| $T$ | set of time periods indexed by $t$ |

**Parameters**

| | |
|---|---|
| $a_{ki}$ | unit transportation cost for raw material between supplier $k$ and candidate site $i$ |
| $bm_{kt}$ | maximum of raw material supply for supplier $k$ in time period $t$ |
| $c_{ij}$ | unit transportation cost for products between candidate site $i$ and customer $j$ |
| $cm_m$ | reference capital cost for manufacture with production capacity level $m$ |
| *d* | range of ambiguity set |
| $d_{jt}$ | demand of customer *j* in time period *t* |
| *f* | density function of probability distribution *P* in ambiguity set *D* |
| $f_0$ | density function of empirical probability distribution $P_0$ |
| *P* | probability distribution in ambiguity set *D* |
| $P_0$ | empirical probability distribution |
| *pex* | unit cost for external product |
| *pi* | unit raw material inventory cost for each time period |
| $pm_m$ | bound of production capacity with production capacity level *m* |
| $vv_i$ | unit production cost in manufacturing facility *i* |
| *α* | worst-case demand fulfillment rate to evaluate the severity of disruptions |
| *β* | conversion rate from raw material to product |
| *η* | raw material degradation rate for each time period |

**Binary variables**

| | |
|---|---|
| $Y_{ir}$ | 0-1 variable. Equal to 1 if a manufacturing facility with production capacity level *r* is built at candidate site *i* |

**Continuous variables**

| | |
|---|---|
| $CA_{ir}$ | production capacity at candidate site *i* with capacity level *r* |
| $CF_i$ | capital cost at candidate site *i* |
| $CO_i$ | production operating cost at candidate site *i* |
| $I_{it}$ | raw material inventory level at candidate site *i* in time period *t* |
| $L_{jt}$ | external product amount at customer *j* in time period *t* |
| $X_{ijt}$ | product transportation flow from candidate site *i* to customer *j* in time period *t* |
| $Z_{kit}$ | raw material transportation flow from supplier *k* to candidate site *i* in time period *t* |

22. Garcia DJ, You F. Supply chain design and optimization: Challenges and opportunities. *Computers & Chemical Engineering.* 2015;81:153-170.
23. Yue D, You F, Snyder SW. Biomass-to-bioenergy and biofuel supply chain optimization: Overview, key issues and challenges. *Computers & Chemical Engineering.* 2014;66:36-56.
24. Verderame PM, Elia JA, Li J, Floudas CA. Planning and scheduling under uncertainty: a review across multiple sectors. *Ind. Eng. Chem. Res.* 2010;49(9):3993-4017.
25. Ning C, You F. Data-Driven Adaptive Robust Unit Commitment Under Wind Power Uncertainty: A Bayesian Nonparametric Approach. *IEEE Transactions on Power Systems.* 2019;34(3):2409-2418.
26. Gong J, You F. Sustainable design and synthesis of energy systems. *Current Opinion in Chemical Engineering.* 2015;10:77-86.
27. Rooney WC, Biegler LT. Optimal process design with model parameter uncertainty and process variability. *AIChE Journal.* 2003;49(2):438-449.
28. Petkov SB, Maranas CD. Design of single-product campaign batch plants under demand uncertainty. *AIChE Journal.* 1998;44(4):896-911.
29. Shang C, You F. A data-driven robust optimization approach to scenario-based stochastic model predictive control. *Journal of Process Control.* 2019;75:24-39.
30. Ben-Tal A, Nemirovski A. Robust optimization - methodology and applications. *Mathematical Programming.* 2002;92(3):453-480.
31. Dupacova J, Growe-Kuska N, Romisch W. Scenario reduction in stochastic programming - an approach using probability metrics. *Math. Program.* 2003;95(3):493-511.
32. Bayraksan G, Love DK. Data-driven stochastic programming using phi-divergences. *The Operations Research Revolution*: INFORMS; 2015:1-19.
33. Ahmed S, Luedtke J, Song YJ, Xie WJ. Nonanticipative duality, relaxations, and formulations for chance-constrained stochastic programs. *Math. Program.* 2017;162(1-2):51-81.
34. Luedtke J, Ahmed S, Nemhauser GL. An integer programming approach for linear programs with probabilistic constraints. *Math. Program.* 2010;122(2):247-272.
35. Jiang RW, Guan YP. Data-driven chance constrained stochastic program. *Math. Program.* 2016;158(1-2):291-327.
36. Yu H, Chung CY, Wong KP, Zhang JH. A chance constrained transmission network expansion planning method with consideration of load and wind farm uncertainties. *IEEE Trans. Power Syst.* 2009;24(3):1568-1576.
37. Cheung WC, Simchi-Levi D. Sampling-based approximation schemes for capacitated stochastic inventory control models, 2015. Located at: Available at SSRN 2615134.
38. Ahmed S, Shapiro A. Solving chance-constrained stochastic programs via sampling and integer programming. *State-of-the-Art Decision-Making Tools in the Information-Intensive Age*: Informs; 2008:261-269.
39. Delage E, Ye YY. Distributionally robust optimization under moment uncertainty with application to data-driven problems. *Oper. Res.* 2010;58(3):595-612.
40. Goh J, Sim M. Distributionally robust optimization and its tractable approximations. *Oper. Res.* 2010;58(4):902-917.
41. Shang C, You F. Distributionally robust optimization for planning and scheduling under uncertainty. *Computers & Chemical Engineering.* 2018;110:53-68.
42. Wiesemann W, Kuhn D, Sim M. Distributionally robust convex optimization. *Oper. Res.* 2014;62(6):1358-1376.36